\magnification=1200
\input amstex
\documentstyle{amsppt}
\parindent=18pt
\topmatter
\title Deformations of weak Fano 3-folds with only terminal
singularities
\endtitle
\author Tatsuhiro Minagawa\endauthor
 \affil Dept. of Mathematical Sciences, University of Tokyo\endaffil
\rightheadtext{Deformations of weak Fano 3-folds}
\abstract
 In this article, we prove that any $\Bbb Q$-factorial weak Fano 3-fold
with only terminal singularities has a smoothing.
\endabstract
\endtopmatter
\head 0. Introduction \endhead 
\definition{Definition 0.1} Let $X$ be a normal Gorenstein projective
variety of dimension 3 over $\Bbb C$ which has only terminal singularities.   
\roster 
 \item If $-K_X$ is ample, we call $X$ a Fano 3-fold.
 \item If $-K_X$ is nef and big, we call $X$ a weak Fano 3-fold.
\endroster
\enddefinition
\definition{Definition 0.2} Let $X$ be a normal Gorenstein projective
variety of dimension 3 with only terminal singularities. Let
$(\varDelta,0)$ be a germ of the 1-parameter unit disk. Let $\frak f
:\frak X \rightarrow (\varDelta,0)$ be a small deformation of $X$ over
$(\varDelta,0)$. We call $\frak f$ a smoothing of $X$ when the fiber
$\frak X_s = \frak f^{-1}(s)$ is smooth for each $s \in
(\varDelta,0)\setminus \{0\}$.
\enddefinition
 We treat the following problem in this paper:
\definition{Problem} Let $X$ be a weak Fano 3-fold with only terminal
singularities.\par
 When $X$ has a smoothing ?  
\enddefinition 
 For the case of Fano 3-fold $X$, $X$ has a smoothing by the result of
Namikawa and Mukai (\cite {Na 3},\cite {Mukai}). Moreover by the method of 
Namikawa, we can show the following theorem:
\proclaim {Theorem 0.3} (Namikawa, Takagi) (Cf. \cite {Na 3},\cite {T})
 Let $X$ be a weak Fano 3-fold with only terminal singularities. Assume
that there exists a birational projective morphism $\pi :X \rightarrow
\bar{X}$ from $X$ to a Fano 3-fold with only canonical singularities
$\bar{X}$ such that $\dim (\pi^{-1}(x))\leq 1$. Then $X$ has a smoothing.
\endproclaim 
\par 
 In this paper, we will show the following theorem:  
\proclaim {Main Theorem} Let X be a weak Fano 3-fold with only terminal
singularities.\roster
 \item The Kuranishi space $Def(X)$ of $X$ is smooth.
 \item There exists $\frak f :\frak X \rightarrow (\varDelta,0)$ a
small deformation of $X$ over $(\varDelta,0)$ such that the fiber $\frak 
X _s =\frak f^{-1}(s)$ has only ordinary double points for 
any $s\in (\varDelta,0)\setminus\{o\}$.
 \item If $X$ is $\Bbb Q$-factorial, then $X$ has a smoothing.
\endroster
\endproclaim
 We remark that if the condition of $(3)$ ``$\Bbb Q$-factorial'' is
dropped, then there is an example that $X$ remains singular under any
small deformation [Cf. (3,7)].
\par
\definition {Acknowledgement} I would like to thank Professor
Y. Namikawa for helpful discussion on the necessity of the condition ``
$\Bbb Q$-factorial''. 
\par
I also express my gratitude to Professor K. Oguiso and H. Takagi  for
helpful conversation.  
\par
I would like to thank Professor Y. Kawamata for useful discussion,
giving me useful suggestions, and encouraging me during the
preparation of this paper.
\enddefinition
\definition{ Notation}\par
 $\Bbb C$ : the complex number field.\par
 $\sim$ : linear equivalence.\par
 $K_X$ : canonical divisor of $X$.\par
 Let $G$ be a group acting on a set $S$. We set $$ S^G:=\{ s\in S \vert
 \,gs=s \text{ for any } g\in G \}.$$\par
 In this paper, $(\varDelta,0)$ means a germ of a 1-parameter unit
disk.\par
 Let $X$ be a compact complex space or a good representative for a
germ, and $\frak g:\frak X \rightarrow (\varDelta,0)$ a 1-parameter small 
deformation of $X$. We denote the fiber $\frak g^{-1}(s)$ for
$s\in (\varDelta,0)$ by $\frak X_s$.\par
 $(\italic {Ens})$ : the category of sets.\par
 Let $\italic k$ be a field. We set $(\italic {Art}_\italic k)$ : the
category of 
Artin local $\italic k$-algebras with residue field $\italic k$.\par
 Let $V$ be a $\Bbb Z$-module. The symbol $V_{\Bbb C}$ means $V
\otimes_{\Bbb Z} \Bbb C$.  
\enddefinition
\head 1. Proof of $(1)$ of Main theorem\endhead
 We use the following theorem of Takagi:
\proclaim  {Theorem 1.1} (Cf. \cite {T})\par
 Let $X$ be a weak Fano 3-fold with only terminal singularities, then
there exists a divisor $D \in \vert -2K_X \vert$ such that $D$ is
smooth.
(In this paper, we call such $D$ a smooth member of $\vert -2K_X\vert$.)
\endproclaim 
\definition {Lemma 1.2}
 Let $X$ be a weak Fano 3-fold with only terminal singularities, $D$ 
a smooth member of $\vert -2K_X \vert$, $A \rightarrow B$ a
surjection in $(\italic {Art}_{\Bbb C})$, $X_B$ an infinitesimal
deformation of $X$ over $B$, and $X_A:=X_B \times _{Spec(B)} Spec (A)$.
Set $D_A \in \vert -2K_{X_A/Spec (A)}$ such that $D_A\vert_X=D$. Then
there exists $D_B \in \vert -2K_{X_B/Spec(B)} \vert$ such that $D_B \vert
_{X_A} = D_A$. 
\enddefinition
\demo{Proof} Since by the Kawamata-Viehweg vanishing theorem we have
 $H^i (X,-2K_X)$ \newline $=0$ for all $i > 0$, we can show above lemma 
as in \cite {Pa, page 63, proof of (iii)}. 
\qed   
\enddemo
\demo{Proof of (1) of Main theorem}
Set $A_n=\Bbb C
[t]/(t^{n+1})$, $\alpha_n:A_{n+1} \rightarrow A_n$, $S_n=Spec (A_n)$.
 Let $D$ be a smooth member of $\vert -2K_X \vert$, and $X_{n+1}$
an infinitesimal deformation of X over $S_{n+1}$, and
$X_n=X_{n+1} \times _{S_{n+1}} S_n$.
By lemma (1.2), there exists $D_{n+1} \in \vert -2K_{X_{n+1}/S_{n+1}} 
\vert$ such that $D_{n+1}\vert _{X} = D$. Set $D_n=D_{n+1} \vert _{X_n}$,
$\pi _{n+1}:Y_{n+1}=Spec (\Cal O _{X_{n+1}} \oplus \Cal O
_{X_{n+1}}(K_{X_{n+1}/S_{n+1}})) \rightarrow X_{n+1}$ a double cover
ramified along $D_{n+1}$, and $\pi _n:Y_n=Spec (\Cal O _{X_n} \oplus \Cal
O_{X_n}(K_{X_n/S_n})) \rightarrow X_n$ a double cover ramified along
$D_n$. We remark that $Y=Y_0$ is a Calabi-Yau 3-fold with only
terminal singularities.
Let $G=\Bbb Z/ 2 \Bbb Z$. We have the following commutative diagram:
$$\CD
  Ext^1_{\Cal O_{Y_{n+1}}}(\Omega^1_{Y_{n+1}/S_{n+1}},\Cal O_{Y_{n+1}})^G
@>T^1_Y(\alpha_n)^G>> Ext^1_{\Cal O_{Y_n}}(\Omega^1_{Y_n/S_n},\Cal
O_{Y_n})^G \\
 @V\beta_{n+1}VV @VV\beta_nV \\
 Ext^1_{\Cal O_{X_{n+1}}}(\Omega^1_{X_{n+1}/S_{n+1}},\Cal O_{X_{n+1}})
@>>T^1_X(\alpha_n)> Ext^1_{\Cal O_{X_n}}(\Omega^1_{X_n/S_n},\Cal
O_{X_n}).
\endCD $$
 We remark that $\beta_{n+1}$ and $\beta_n$ are defined because
$\pi_{n+1}$ and $\pi_n$ are finite morphisms. (Cf. \cite {Na  1, 
Proposition 4.1}). By \cite {Na 1, proof of Theorem 1, page 431}, 
$T^1_Y(\alpha_n)$
is surjective. Thus $T^1_Y(\alpha_n)^G$ is a surjection
because $G$ is finite. $\beta_n$ is a surjection by
lemma (1.2) and we have that $T^1_X(\alpha_n)$ is also surjective. By
$\Bbb T ^1$-lifting criterion (Cf \cite {Ka 1}, \cite {Ka 2}),
we proved (1) of Main theorem.
\qed
\enddemo 
\head 2. Proof of (2) of Main theorem \endhead
 We use the result of Namikawa and Steenbrink on deformations of 
Calabi-Yau 3-folds to prove (2) of Main theorem.
Let $Y$ be a Calabi-Yau 3-fold with only terminal singularities,
$\{q_1, q_2, \dots, q_n\}=Sing(Y)$,
$\nu : \Tilde{Y} \rightarrow Y $ be a good resolution of $Y$, and 
$E_i= \nu ^{-1}(q_i)$. (``good'' means the restriction of 
$\nu : \nu ^{-1}(V) \rightarrow V$ is an isomorphism and
its exceptional divisor $E_i$ is simple normal crossings for each $i$.)
\proclaim {Proposition 2.1}(Cf. \cite {Na-St}  )
 If $(Y,q_i)$ is not the ordinary double point, then the homomorphism 
$\iota_i:H^2_{E_i}(\Tilde{Y},\Omega ^2_{\Tilde{Y}}) \rightarrow
H^2(\Tilde{Y}, \Omega ^2_{\Tilde{Y}})$ is not injective.
\endproclaim
\demo{Proof of (2) of Main theorem}
 By theorem (1.1), there exists a smooth member $D$ of $\vert -2K_X
\vert$. We remark that $D \cap Sing(X)=\emptyset$.
Let $\{p_1, p_2, \dots, p_n \}=Sing(X)$, and $\pi:Y=Spec
(\Cal O _X \oplus \Cal O _X(K_X)) \rightarrow X$ be a double cover
ramified along $D$. Then $Y$ is a Calabi-Yau 3-fold with only
terminal singularities. Let $G=\Bbb Z / 2 \Bbb Z =\{id_X, \sigma\}$,
$\pi ^{-1}(p_i)=\{q_{i1},q_{i2}\}$.
Then we have that $Sing(Y)=\{q_{ij}\vert i=1, 2, \dots,
n, j=1, 2\}$ because $D$ is smooth. Let $Y_{ij}$ be a sufficiently
small open neighborhood of $q_{ij}$, $V_{ij}=Y_{ij}\setminus
\{q_{ij}\}$, $U=X \setminus Sing(X)$ , and $V=Y \setminus Sing(Y)$.
Let $\nu:\Tilde {Y}
\rightarrow Y$ be a $G$-equivariant good resolution of $Y$, and
$E_{ij}=\pi^{-1}(q_{ij})$. Let $\omega \in H^0(\omega _Y)$ be a nowhere 
vanishing section. We remark that $\sigma(\omega)=-\omega$.
We consider the following commutative diagram: 
$$\CD
  [H^1(V,\Omega^2 _V)]^{[-1]} @>\alpha^{\prime}>> 
[\oplus_{i,j}H^2_{E_{ij}}(\Tilde
{Y},\Omega^2_{\Tilde {Y}})]^{[-1]} @>\iota>> 
[H^2(\Tilde{Y},\Omega^2_{\Tilde {Y}})]^{[-1]}\\
  @AA\wr \,\tau A @AA\oplus_{i,j} \tau _{ij} A \\
  H^1(V,\Theta _V)^G @>\alpha>> \oplus _{i,j} H^1(V_{ij},\Theta _{V_{ij}})^G
\endCD $$
where $F^{[-1]}=\{x\in F \vert \sigma(x)=-x \}$	for a $\Bbb C$-vector
space $F$ with a $G$-action.\par
 By proposition (2.1), $\iota_{ij}$ is not injective if $(X,p_i) \simeq
(Y,q_{ij})$ is not the ordinary double point. So there exists an element
$\eta ^{\prime}\in [H^1(V,\Theta_{V})]^{[-1]}$ such that
$\alpha^{\prime}(\eta)_{ij} \not= 0$ for any $i,j$ 
where $(X,p_i) \simeq (Y, q_{ij})$ is not the
ordinary double point. Let $\eta \in H^1(V,\Theta_{V})^G$ such that
$\tau(\eta) = \eta^{\prime}$. Let $\beta : [H^1(V,\Theta_V)]^G
=Ext^1_{\Cal O_Y}(\Omega ^1_Y, \Cal O_Y)^G \rightarrow
Ext^1_{\Cal O _X}(\Omega^1_X, \Cal O_X) = H^1(U,\Theta_U)$ 
be the homomorphism defined in the proof of (1) of Main 
theorem. By (1) of Main theorem, there exists a small deformation of $X$
over $(\varDelta,0)$ $\frak f : \frak X \rightarrow (\varDelta, 0)$
which is a realization of $\beta(\eta)$. 
Using the method of Namikawa
(Cf. \cite{Na 2 , Theorem 5}, \cite{Na-St, Theorem (2.4)}),
we can reach a smooth 3-fold by small
deformations by continuing the process above.
\qed
\enddemo
\definition {Definition 2.2} Let $X$ be a normal $\Bbb Q$-Gorenstein
projective variety of dimension 3 over $\Bbb C$ which has only terminal
singularities.
\roster
 \item The index $i_p$ of
a singular point $p \in X$ is defined by $$ i_p:=min\{ m\in\Bbb N \vert
\, mK_X \text{ is a Cartier divisor near $p$} \}.$$
 \item The sigular index $i(X)$ of $X$ is defined by $$ i(X):
=min\{ m\in\Bbb N \vert \, mK_X \text{ is a Cartier divisor} \}.$$
 \item If $-K_X$ is ample, we call $X$ a $\Bbb Q$-Fano 3-fold.
 \item If $-K_X$ is nef and big, we call $X$ a weak $\Bbb Q$-Fano
3-fold.
\endroster
\enddefinition
 We considered deformations of $\Bbb Q$-Fano 3-folds in \cite {Mi}.
The method of (2) of Main theorem is also useful for weak $\Bbb Q$-Fano
3-folds of singular index 2.
\definition {Definition 2.3}  Let $(X,p)$ be a germ of a 3-dimensional
terminal singularity and $G=\Bbb Z / 2 \Bbb Z$. 
We call $(X,p)$ a quotient of the ordinary double
point if $(X,p)$ is isomorphic to the singularity of the following type:
 Let $x_1,$ $x_2,$ $x_3,$ $x_4$ be coordinates of the germ $(\Bbb C
^4,0)$. We define a $G$-action on $(\Bbb C^4,0)$ by $
x_1 \mapsto -x_1$, $x_2 \mapsto -x_2$, $x_3 \mapsto x_3$, $x_4 \mapsto
-x_4$. $(X,p) \simeq \{x_1^2 + x_2^2 + x_3^2 + x_4^2 =0 \vert (\Bbb C
^4,0)\}/G$.
\enddefinition
\proclaim {Theorem 2.4}
 Let $X$ be a weak $\Bbb Q$-Fano 3-fold with only terminal
singularities of singular index $i(X)=2$, and assume that there
exists a smooth member of 
$\vert -2K_X \vert$. Then there exists a small deformation of $X$ over
$(\varDelta,0)$ $\frak f:\frak X \rightarrow (\varDelta,0)$ such that the 
fiber $\frak X _s =\frak f^{-1}(s)$ has only ordinary double points or
quotients of ordinary double points for any $s \in (\varDelta,0)
\setminus \{0\}$.
\endproclaim
 To prove this theorem, we use an analogous proposition of 
proposition (2.1).
 Let $X$ be a weak $\Bbb Q$-Fano 3-fold with only terminal
singularities of singular index $i(X)=2$, and assume that there
exists a smooth member of 
$\vert -2K_X \vert$. Let $D$ be a smooth member of $\vert -2K_X \vert$.
Let $\pi:Y=Spec (\Cal O _X \oplus \Cal O _X(K_X)) \rightarrow X$ 
be a double cover ramified along $D$.
Then $Y$ is a Calabi-Yau 3-fold with only terminal singularities.
Let $p \in X$ be a singularity of index $i_p =2$, and $\pi^{-1}(p)={q}$.
We remark that $\pi \vert _{(Y,q)}: (Y,q) \rightarrow (X,p)$ is a
canonical cover of $(X,p)$. Let $G=\Bbb Z/ 2 \Bbb Z$. $\nu : 
\Tilde{Y} \rightarrow Y $ be a $G$-equivariant good resolution of $Y$,
and $E= \nu ^{-1}(q)$. We know the following proposition
which is analogous to proposition (2.1) and is a result of Namikawa.
\definition {Proposition 2.5} (Cf. \cite{Na 4})
If $(X,p)$ is a singular point of index $i_p=2$, and if $(Y,q)$ is not
the ordinary double point, then the homomorphism 
$\iota^{[-1]}:H^2_{E}(\Tilde{Y},\Omega ^2_{\Tilde{Y}})^{[-1]} 
\rightarrow H^2(\Tilde{Y}, \Omega ^2_{\Tilde{Y}})^{[-1]}$ is 
not injective.
\enddefinition
 This proposition leads us to theorem (2.4) by the same method of the
proof of (2) of Main theorem.    				
\head 3. Proof of (3) of Main theorem \endhead
 We first prove the following theorem to prove (3) of Main theorem.
\proclaim {Theorem 3.1} 
 Let $X$ be a weak Fano 3-fold with only terminal singularities. Assume
that $X$ is $\Bbb Q$-factorial, then there exists a divisor $S \in \vert
-K_X \vert$ such that $S$ is smooth.
\endproclaim
 To prove theorem (3.1), we use some known results as follows.
\definition {Definition 3.2}
 Let $X$ be a weak Fano 3-fold with only terminal singularities. Fano
index of $X$ is defined by
 $$F(X)=max\{ r \in \Bbb N \vert \text{ there exists a Cartier divisor }
H \text{ such that } -K_X \sim rH \}.$$
\enddefinition    
\proclaim {Theorem 3.3} (Reid, Shin) (Cf. \cite{Re 2}, \cite{Shi})
 Let $X$ be a Fano 3-fold with only canonical singularities. Then we
have,
\roster 
 \item $\dim Bs \vert -K_X \vert \leq 1$, 
 \item if $F(X)>1$ then $Bs \vert -K_X \vert= \emptyset$,
 \item if $\dim \vert -K_X \vert=1$ then a general member of $\vert -K_X
\vert$ is smooth at base points of $\vert -K_X \vert$, and
 \item if $\dim Bs \vert -K_X \vert =0$ then $Bs \vert -K_X \vert=\{p\}$
one point, a general member of $\vert -K_X \vert$ has the ordinary double
point at $p$, and $p \in Sing(X)$.
\endroster
\endproclaim
\proclaim {Theorem 3.4} (Mella) (Cf. \cite {Me, Theorem (2.4)})
  In the case of (1) of theorem (3.3), if $p \in Sing(X)$ is a terminal
singularity, then $X \cong X_{2,6} \subset \Bbb
P(1,1,1,1,2,3)$. Moreover for any Zariski open set $U$ containing $p$,
$U$ is not $\Bbb Q$-factorial.
\endproclaim
\proclaim {Theorem 3.5} (Reid, Ambro) (Cf. \cite {Re 2}\cite {Am})
 Let $X$ be a weak Fano 3-fold with only canonical singularities, then
a general member of $\vert -K_X \vert$ has only canonical singularities.
\endproclaim
\demo {Proof of theorem (3.1)}    
 Let $\pi : X \rightarrow \bar{X}$ be a multi-anti-canonical morphism,
then $\bar{X}$ is a Fano 3-fold with only canonical singularities, and
$\pi$ is crepant ($K_X=\pi^{\ast}(K_{\bar{X}})$).\par 
 In the case of $Bs\vert -K_{\bar{X}} \vert =\emptyset$ or 
$\dim Bs \vert -K_{\bar{X}} \vert=1$,
then a general member of $\vert -K_{\bar{X}} \vert$ is smooth at its
base point by theorem (3.3), and there exists a divisor
$S \in \vert -K_X \vert$ such that $S$ is smooth by theorem (3.5).\par
 In the case of $\dim Bs\vert -K_{\bar{X}} \vert=0$ (in this  case $Bs\vert
-K_X\vert = \{p\}$ by theorem (3.3.4)), there exists a divisor $\bar{S}
\in \vert -K_{\bar{X}} \vert$ which has the ordinary double point at $p$ 
such that $S=\pi^{\ast}(\bar{S})$ has only canonical singularities. 
If we can not take a smooth $S$, then $\pi \vert _S : S \rightarrow
\bar{S}$ is an isomorphism  near $p$ because $p$ is the ordinary double
point. Then there exists a Zariski open set $U$ containing $\pi ^{-1}(p)$
such that $\pi \vert _U : U \rightarrow \bar{X}$ is an open immersion.
So $p \in \bar{X}$ is terminal. By theorem (3.4), $\pi(U)$ is not $\Bbb
Q$-factorial. Thus $U$ is not $\Bbb Q$-factorial and $X$ is not $\Bbb
Q$-factorial which is a contradiction.
\qed
\enddemo
 By (2) of Main theorem, the following theorem  is enough
to prove (3) of Main theorem.
\proclaim {Theorem 3.6} Let $X$ be a weak Fano 3-fold with only ordinary 
double points. Assume that $X$ is $\Bbb Q$-factorial. Then $X$ has a
smoothing. 
\endproclaim
\demo{Proof}
 Let $\nu:\Tilde{X} \rightarrow X$ be a small resolution of $X$,
$\{p_1,p_2,\dots,p_n\}=Sing(X)$, $U=X-Sing(X)$, $X_i$ a sufficiently
small open neighborhood of $p_i$, $U_i=X_i \setminus \{p_i\}$, 
and $C_i= \nu^{-1}(p_i)$. Since
$H^1(X,\Omega^1_X) \simeq H^1(X,\nu_{*} \Omega^1_{\Tilde{X}}) $ 
(Cf. \cite{Na 1, Lemma (2.2)}), We have the following commutative 
diagram of exact sequences:
$$\CD
  0 @>>> H^1(X,\Omega^1_X)@>\alpha _1>>
H^1(\Tilde{X},\Omega^1_{\Tilde{X}})@>\alpha_2>>
H^0(X,R^1\nu_{*}\Omega^1_{\Tilde{X}})\\
  @. @AA \lambda_1A @AA \lambda_2A @AA \lambda_3A \\  
  0 @>>> H^1(X, \Cal O ^*_X)_{\Bbb C} @>\beta_1>>
H^1(\Tilde{X}, \Cal O ^*_{\Tilde{X}})_{\Bbb C} @>\beta_2>>
H^0(X,R^1\nu_*\Cal O^*_{\Tilde{X}})_{\Bbb C} 
\endCD $$ 
$\lambda_2$ is surjective because $h^2(\Tilde{X},\Cal O _{\Tilde{X}})=0$, 
and $\beta_1$ is also surjective because $X$ is $\Bbb Q$-factorial and
$\nu$ is small. Thus we have that $\alpha_2$ is the zero map, and its
dual $\oplus_{i=1}^n H^2_{C_i}(\Tilde{X}, \Omega^2_{\Tilde{X}}) \rightarrow
H^2(\Tilde{X}, \Omega^2_{\Tilde{X}})$ is also the zero map. By theorem
(3.1), there exists $D \in \vert -K_X \vert$ a smooth member of $\vert
-K_X \vert$. Then $D \cap Sing(X)=\phi$. We consider the following
commutative diagram defined by $\nu^* D$:
$$\CD
 \oplus_{i=1}^n H^2_{C_i}(\Tilde{X},\Omega^2_{\Tilde{X}}) @>\oplus_i
\delta_i >> \oplus _{i=1} ^n H ^2 _{C_i} (\Tilde{X},\Theta _{\Tilde{X}}) \\
 @VVV  @VV{\oplus_i \iota_i}V \\
 H^2(\Tilde{X},\Omega^2_{\Tilde{X}}) @>>>
H^2(\Tilde{X},\Theta_{\Tilde{X}}).
\endCD $$
 $\delta_i$ is an isomorphism for any $i$, and we have that $\iota_i$ is 
the zero map for any $i$. We consider the following exact commutative 
diagram:
$$\CD
 H^1(U,\Theta_U) @>\gamma^{\prime}>> \oplus _{i=1}^n
H ^2 _{C_i}(\Tilde{X},\Theta _{\Tilde{X}}) @>\oplus_i
\iota_i>>H^2(\Tilde{X},\Theta_{\Tilde{X}}) \\
@| @AAA \\
 H^1(U,\Theta_U) @>>\gamma> \oplus_{i=0}^n H^1(U_i,\Theta_{U_i}).
\endCD $$ 
 Then there exists an element $\eta \in H^1(U,\Theta_U)$ such that
$\gamma^{\prime}(\eta)_i \not= 0$ for any $i=1,2,\dots,n$. Thus
$\gamma(\eta)_i \not= 0$ for any $i=1,2,\dots,n$. By (1) of Main theorem,
there exists a small deformation of $X$ over $(\varDelta,0)$ $\frak f:
\frak X \rightarrow (\varDelta,0)$ which is a realization of
$\eta$. Then $\frak f$ is a smoothing of $X$.
\qed
\enddemo     
\definition {Example 3.7}
 Let $\bar{X}$ be the projective cone over the smooth del Pezzo surface $S$ of
degree 8. Then $\bar{X}$ is a Gorenstein Fano 3-fold with $\rho = 1$
which has only  one Gorenstein rational singularity $\bar{p}$ at its vertex.
Let $f: Z \rightarrow \bar{X}$ be the blowing-up at $\bar{p}$, then $f$ is a
crepant resolution of $\bar{X}$ and $Z \simeq Proj(\Cal O _S \oplus
\omega ^{-1}_S)$. Let $E$ be an exceptional divisor of $f$ which is
isomorphic to $\Bbb F _1$, and $C$ be the $(-1)$-curve on $E$. Then $Z$ is a
weak Fano 3-fold with $(-1,-1)$-curve $C$. Let $\nu: Z \rightarrow X$ be
a birational contraction which contracts $C$. Then $X$ is a weak Fano
3-fold which has only one ordinary double point $\nu (C) =p$. Let
$F=\nu(E)$, then $F \simeq \Bbb P ^2$ passing through $p$. So $X$ is not 
$\Bbb Q$-factorial. We have that $X$ is not smoothable, in fact there exists a
sufficiently small open neighborhood $U$ of $F$ $(F \subset U)$ which is not
smoothable by \cite{Na 5, Proposition (1.3)}.
\enddefinition

\enddocument